\documentclass[10pt]{article}

\usepackage{MacrosRM}
\newcommand{\dd}{\mathrm{d}}
\newcommand\numberthish{\addtocounter{equation}{1}\tag{\theequation}}
    
\begin{document}

\title{\textbf{Conditional Backward Propagation of Chaos}}
\author{\small\textsc{Rémi Moreau} \thanks{Univ. Rennes, IRMAR, 35000, Rennes, \textsc{France}. remi.moreau@ens-rennes.fr}}
\date{}
\maketitle

\begin{abstract}
In this paper, we first investigate the well-posedness of a backward stochastic differential equation where the driver depends on the law of the solution conditioned to a common noise. Under standard assumptions, we show that existence and uniqueness, as well as integrability results, still hold. We also study the associated interacting particles system, for which we prove propagation of chaos, with quantitative estimates on the rate of convergence in Wasserstein distance.
\end{abstract}

\bigskip

\textbf{MSC-classification}: 60F25, 60H10, 60H20 
\bigskip

\section{Introduction}

The theory of backward stochastic differential equations (BSDEs) was first introduced by Pardoux and Peng, among others. In \cite{PP1990}, they solve the following problem, where $W$ denotes a Brownian motion on a given probability space equipped with an adapted filtration :
\begin{equation}\label{EqIntro}
    Y_t = \xi + \int_t^T f_s(Y_s,Z_s)\, \dd s - \int_t^T Z_s \, \dd W_s
\end{equation}
for a driver $f$ which is Lipschitz continuous, and a square integrable terminal condition $\xi$. To solve \eqref{EqIntro} means to prove existence and uniqueness of an adapted process $(Y,Z)$, with fitting integrability properties, satisfying \eqref{EqIntro}. 
\medskip

The last few decades have seen many extensions, some with applications in stochastic control or mathematical finance. The hypothesis on the driver $f$ can be weakened (see for example \cite{briandpolynomial} where the driver $f$ is polynomial in the $y$ component, or \cite{Kob97} where it is quadratic in the $z$ component), the solution follows the integrability of the terminal condition (see \cite{HuDelyon03} for existence and uniqueness when $\xi\in L^p$ for instance), constraints can be applied to the solution, which gives rise to reflected BSDEs and comparison results (see \cite{ElK97}, \cite{Ham00}), or even to the law of the solution (\cite{BCCdRH}, \cite{BriandElieHu}); BSDEs with jumps can also be considered (\cite{Roy06}), to cite a few of these various existing developments.
\medskip

On the other hand, the theory of propagation of chaos, led by the work of Kac, aims at studying the behavior of a group of interacting individuals or particles when the size of the group goes to infinity. The initial independance is "propagated" in the sense that the particles become indepedant at the limit, the interaction becoming weaker as the size of the group grows. Notes by Sznitman (\cite{Sznitman}) and a more recent paper by Fournier and Guillin (\cite{FournierGuillin}) give a refinement to the weak convergence given by the Glivenko-Cantelli theorem for measures admitting high-enough moments, namely in terms of convergence in Wasserstein spaces.
\medskip

The result with common noise in the forward SDE setting can be found in different contexts, for example in \cite[Theorem 2.12]{CD2} for a quite general case, or in \cite{BCCdRH} where the system of interacting particles is furthermore subject to a constraint in law. 
\medskip

In the backward setting, a simplified, without environmental noise, version of our main result is given in \cite{LauriereTangpi} in the backward framework. The terminal conditions are now fixed, independant, and the interaction is modeled by BSDEs through the empirical distribution.
\medskip

Our aim in this paper is to extend the results to the following corresponding problem :
\begin{equation}\label{PbInit}
    \dd Y_t = - f_t(Y_t,Z_t,Z_t^0,\mu_t) \, \dd t + Z_t\, \dd W_t + Z_t^0 \, \dd W_t^0,~~ Y_T = \xi
\end{equation}
where $(\mu_t)$ denotes the conditional law of the solution $Y$ with respect to the common noise $W^0$. We raise awareness on the fact that, without this dependency in the conditional law, the result is simply a direct application of the works of Lauriere and Tangpi. As in \cite{LauriereTangpi}, the result relies heavily on the coupling technique, classical BSDE estimates and well-known conditions (see for example \cite{CD1}) allowing to get rates of convergence for the propagation of chaos.
\medskip

The paper is organized as follows : we start by proving that the problem $\eqref{PbInit}$ is well-posed and admits a unique solution under suitable conditions. We then make use of classical adapted results on the moments of the solution to prove two results of backward propagation of chaos for the empirical measure of the conditioned system, under two different sets of hypotheses. The last part is devoted to the proof of pointwise and uniform convergence for the corresponding system of interacting particles.
\bigskip

\section{Setting and first results}
\noindent Throughout this paper, we fix $T>0$ a fixed finite time horizon and we denote by $| \cdot |_E$ the euclidian norm on the space $E$ (we drop the subscript when the space is $\R^n$). Let us consider two complete probability spaces $(\Omega^0, \FF^0, \P^0)$ and $(\Omega^1, \FF^1, \P^1)$, endowed with $\F^0 = (\FF_t^0)_{t\geq 0}$ and $\F^1 = (\FF_t^1)_{t\geq 0}$ respectively, two right-continuous complete filtrations. As in \cite{CD2}, we consider the complete product space $(\Omega, \FF, \P)$ where $\Omega = \Omega^0\times \Omega^1$, and $(\FF, \P)$ is the completion of $(\FF^0\otimes\FF^1, \P^0\otimes \P^1)$, endowed with $\F = (\FF_t)_{t\geq 0}$ the right-continuous complete augmentation of $(\FF_t^0\otimes \FF_t^1)$, on which one can define $\LL^1(X)$ for $X$ a random variable on the product space. As shown in \cite[Lemma 2.4]{CD2}, the law $\LL^1(X)$ can be seen as the conditional law with respect to $\F^0$.
Let us consider the following problem :
\begin{equation}\label{eq1}
    Y_t = \xi + \int_t^T f(s,Y_s,Z_s,Z_s^0,\mu_s)\, \dd s - \int_t^T Z_s\, \dd W_s - \int_t^T Z_s^0\, \dd W_s^0,
\end{equation}
where $W$ and $W^0$ are independent Brownian motions on $\Omega^1$ and $\Omega^0$ respectively, and $(\mu_s)_{0\leq s\leq T} = (\mathcal{L}^1(Y_s))$. We introduce the following notations, for $p\geq 1$, which will be useful in the rest of the paper.

\begin{itemize}
    \item $L^p$ is the collection of $\R^m$-valued $\FF_T$-measurable variables $\xi$ satisying 
    \[\left\| \xi \right\|_p^p = \E \left[ |\xi|^p \right] < \infty;\]
    \item $\HH^p$ is the collection of $\R^{m\times d}$-valued $\F$-progressively measurable processes $z$ satisfying 
    \[\left\| z \right\|_{\HH^p}^p = \E \left[ \left( \int_0^T |z_s|^2\, \dd s\right)^{\frac{p}{2}} \right] < \infty ;\]
    \item $\SS^p$ is the collection of $\R^{m}$-valued $\F$-progressively measurable processes $y$ satisfying 
    \[\left\| y \right\|_{\SS^p}^p = \E \left[ \sup\limits_{t\leq T} |y_t|^p \right] < \infty .\]
\end{itemize}
We also denote $\E^0$ the conditional expectation with respect to the common noise $W^0$. In what follows, we make use of the following conditions on the terminal condition $\xi$ and on the driver $f$.
\begin{enumerate}
    \item[\textbf{(H1)}] The terminal condition $\xi$ is in $L^2$.
    \item[\textbf{(H2)}] We have the following usual Lipschitz condition and the process $(f(t,0,0,0,\delta_0))$ belongs to $\mathcal{H}^2$ :
    \begin{align*}
    \forall (t,y,z,\tilde{z},\mu), (t,y',z',\tilde{z}',\nu),~ |f&(t,y,z,\tilde{z},\mu) - f(t,y',z',\tilde{z}',\nu)|\\ &\leq C_f \left( | y-y'| + | z-z'| + |\tilde{z}-\tilde{z}'| + W_2(\mu,\nu)\right) ~\P-a.s.\\
    \E\left[ \int_0^T |f(s,0,0,0,\delta_0)|^2\, \dd s\right] < +\infty
    \end{align*}
\end{enumerate}

\begin{Def}
    A solution to \eqref{eq1} is a triple of progressively measurable processes $(Y,Z,Z^0)$ taking values in $\R^m \times \R^{m\times d} \times \R^{m\times d}$ satisfying \eqref{eq1} with $Y\in \SS^2$, $Z,Z^0\in \HH^2$.
\end{Def}

\begin{Thm}
    If (H1) and (H2) hold, then the equation \eqref{eq1} admits a unique solution $(Y,Z,Z^0)\in \SS^2 \times \HH^2 \times \HH^2$.
\end{Thm}

\begin{Rmq}
    This is a generalization of the case where there is no common noise, in which $(\mu_t)$ is replaced by the law of the solution on the whole probability space.
\end{Rmq}

\begin{proof}
    We recall the sketch of this proof for the reader's convenience, with extra emphasis on what differs from the classical case. Consider $\Phi$ the mapping $(U,V,V^0) \mapsto (Y,Z,Z^0)$ where $(Y,Z,Z_0)$ is the unique solution to the following classical BSDE (with $(W,W^0)$ seen as a $2d$-dimensional Brownian motion) :
    \[\dd Y_t = -f(t,U_t,V_t,V_t^0,\LL^1 (U_t)) \, \dd t + Z_t\, \dd W_t + Z_t^0\, \dd W_t^0,~~ Y_T = \xi.\]
    Fix $(U^1,V^1,V^{0,1})$ and $(U^2,V^2,V^{0,2})$ and denote $(Y^i,Z^i,Z^{0,i}) = \Phi(U^i,V^i,V^{0,i})$ for $i\in\{1,2\}$. In the following, we write $\Delta f_t = f(t,U_t^1,V_t^1,V_t^{0,1},\LL^1 (U_t^1)) - f(t,U_t^2,V_t^2,V_t^{0,2},\LL^1 (U_t^2))$ and $\Delta M_t = M_t^1 - M_t^2$, for $M=U^i,V^i,V^{0,i},Y^i,Z^i,Z^{0,i}$. Let us start by expanding $e^{\alpha t} |\Delta {Y}_t|^2$ using the Itô formula :
    \begin{align*} e^{\alpha t} |\Delta Y_t|^2 + &\int_t^T e^{\alpha s} (\alpha |\Delta Y_s|^2 + |\Delta Z_s|^2 + |\Delta Z_s^0|^2) \, \dd s \\ &= 2\int_t^T e^{\alpha s} \Delta Y_s \cdot \left( \Delta f_s \, \dd s - \Delta Z_s \, \dd W_s - \Delta Z_s^0\, \dd W_s^0\right). \end{align*}
    Now, taking the conditional expectation with respect to $\FF_t$, we get :
    \[e^{\alpha t} |\Delta Y_t|^2 + \E \left[ \int_t^T e^{\alpha s} (\alpha |\Delta Y_s|^2 + |\Delta Z_s|^2 + |\Delta Z_s^0|^2) \, \dd s \mid \FF_t\right] = 2\E \left[ \int_t^T e^{\alpha s} \Delta Y_s \cdot \Delta f_s\, \dd s \mid \FF_t \right].\]
    For any $\varepsilon > 0$, by exploiting the classical inequality $ab \leq a^2 \varepsilon + b^2 \varepsilon^{-1}$ and recalling the Lipschitz condition on $f$,
    \begin{align*}
        &e^{\alpha t} |\Delta Y_t|^2 + \E \left[ \int_t^T e^{\alpha s} (\alpha |\Delta Y_s|^2 + |\Delta Z_s|^2 + |\Delta Z_s^0|^2) \, \dd s \mid \FF_t\right] \\ &\leq \frac{1}{\varepsilon} \E \left[ \int_t^T e^{\alpha s} |\Delta Y_s|^2\, \dd s \mid \FF_t\right] \\ & \hspace*{2cm}+ 8\varepsilon C_f^2 \E \left[ \int_t^T e^{\alpha s} \left( |\Delta U_s|^2 + |\Delta V_s|^2 + |\Delta V_s^0|^2 + W_2^2(\LL^1(U_t^1), \LL^1(U_t^2))\right)\, \dd s \mid \FF_t\right].
    \end{align*}
    Taking the expectation (with respect to $\P$), we get :
    \begin{align*}
        \E &\left[ \int_t^T e^{\alpha s} \left( \left(\alpha - \frac{1}{\varepsilon} \right) |\Delta Y_s|^2 + |\Delta Z_s|^2 + |\Delta Z_s^0|^2 \right) \, \dd s \right] \\ &\leq 16 C_f^2 \varepsilon \E \left[ \int_t^T e^{\alpha s} |\Delta Y_s|^2\, \dd s \right] + 8\varepsilon C_f^2 \E \left[ \int_t^T e^{\alpha s} \left( |\Delta U_s|^2 + |\Delta V_s|^2 + |\Delta V_s^0|^2 \right)\, \dd s\right],
    \end{align*}
    where we used the following inequality, stemming from the definition of Wasserstein distances :
    \begin{equation} \E\left[ W_2^2(\LL^1(U_t^1), \LL^1(U_t^2)) \right] \leq \E\left[\E^0 \left[ |U_t^1 - U_t^2|^2 \right] \right] = \E\left[ |U_t^1 - U_t^2|^2 \right].\end{equation}
    By choosing $\varepsilon$ and $\alpha$ such that $16C_f^2 \varepsilon < 1$ and $\alpha\geq 1+\varepsilon^{-1}$, there exists $C < 1$ such that :
    \[\E \left[ \int_t^T e^{\alpha s} \left( |\Delta Y_s|^2 + |\Delta Z_s|^2 + |\Delta Z_s^0|^2 \right) \, \dd s \right] \leq C \E \left[ \int_t^T e^{\alpha s} \left( |\Delta U_s|^2 + |\Delta V_s|^2 + |\Delta V_s^0|^2 \right)\, \dd s\right]\]
    As a consequence, the mapping $\Phi$ defined above admits a unique fixed point, which solves (uniquely) the equation $\eqref{eq1}$.
\end{proof}

\begin{Rmq}
    The result of existence and uniqueness above also gives the following uniform bound :
    \[\E\left[ \sup\limits_{0\leq t\leq T} |Y_t|^2 \right] < +\infty.\]
\end{Rmq}

\noindent Let us now turn to a few useful inequalities on the moments of the solution.

\begin{Prop}[] \label{PropLp}
    For $p\geq 2$, assume furthermore that $\xi\in L^p$ and that $(f(s,0,0,0,\delta_0))\in \HH^p$, which means that the following holds :
    \[\E\left[ \int_0^T |f(s,0,0,0,\delta_0)|^p\, \dd s\right] < +\infty.\]
    Then, the solution $Y$ belongs to $\SS^p$, that is :
    \[\E\left[ \sup\limits_{t\leq T} |Y_t|^p\right] < +\infty.\]
\end{Prop}

\begin{proof}
    Let us start by applying Itô's formula to $e^{\alpha t} |Y_t|^p$, with $p\geq 2$ and $\alpha >0$, which gives :
    \begin{align*}
        e^{\alpha t} |Y_t|^p &+ \frac{p(p-1)}{2} \int_t^T e^{\alpha s} |Y_s|^{p-2} \left(|Z_s|^2 + |Z_s^0|^2\right) \dd s \\ &= |\xi|^p + p\int_t^T e^{\alpha s} |Y_s|^{p-2} Y_s \cdot f(s,Y_s,Z_s,Z_s^0,\mu_s)\, \dd s \\
        & \hspace*{1cm}- p \int_t^T  e^{\alpha s} |Y_s|^{p-2} Y_s\cdot \left( Z_s\, \dd W_s + Z_s^0\, \dd W_s^0 \right) - \alpha \int_t^T e^{\alpha s}|Y_s|^p\, \dd s.
    \end{align*}
    The Lipschitz continuity of $f$ and Young's inequality provide the following inequalities :
    \begin{align*}
        \int_t^T e^{\alpha s} &|Y_s|^{p-2} Y_s \cdot f(s,Y_s,Z_s,Z_s^0,\mu_s)\, \dd s \\ & \leq \int_t^T e^{\alpha s} |Y_s|^{p-2} |Y_s| \left( |f(s,0,0,0,\delta_0)| + |Y_s| + |Z_s| + |Z_s^0| + W_2(\mu_s,\delta_0) \right) \dd s \\
        &\leq C_p \int_t^T e^{\alpha s} \left[ |f(s,0,0,0,\delta_0)|^p + \left( 1+ \frac{1}{\varepsilon} \right) |Y_s|^p + \varepsilon |Y_s|^{p-2} \left(|Z_s|^2 + |Z_s^0|^2\right) + \varepsilon W_p^p(\mu_s,\delta_0) \right]\dd s.
    \end{align*}
    since $W_2 \leq W_p$. Coming back to the first identity, for $\varepsilon$ such that $4C_p\varepsilon \leq p(p-1)$ and $\alpha$ big enough, we get :
    \begin{align*}
        & e^{\alpha t} |Y_t|^p + \frac{p(p-1)}{4} \int_t^T e^{\alpha s} |Y_s|^{p-2} \left(|Z_s|^2 + |Z_s^0|^2\right) \dd s \\ &\hspace{2cm} \leq C_{p,T} \left( |\xi|^p + \int_t^T \left( |f(s,0,0,0,\delta_0)|^p + \varepsilon W_p^p(\mu_s,\delta_0) \right)\dd s\right) \numberthish \label{eqn2}\\
        &\hspace*{4cm} - p \int_t^T  e^{\alpha s} |Y_s|^{p-2} Y_s\cdot \left( Z_s\, \dd W_s + Z_s^0\, \dd W_s^0 \right)
    \end{align*}
    Using the following property of the Wasserstein distance $\E\left[ W_p^p(\mu_s, \delta_0) \right] \leq \E\left[\E^0 \left[ |Y_t - 0|^p \right] \right] = \E\left[ |Y_t|^p \right]$ previously mentioned and Gronwall's inequality, the above inequality becomes, taking the expectation first and the supremum in time second :
    \begin{align*} &\sup\limits_{t\leq T} \E\left[ e^{\alpha t} |Y_t|^p + \frac{p(p-1)}{4} \int_t^T e^{\alpha s} |Y_s|^{p-2} \left(|Z_s|^2 + |Z_s^0|^2\right) \dd s\right] \\ &\hspace*{3cm}\leq C_{p,T} \E\left[ |\xi|^p + \int_0^T |f(s,0,0,0,\delta_0)|^p\, \dd s \right].\end{align*}
    In order to swap the supremum and the expectation, we go back to the equation \eqref{eqn2} and apply the Gronwall's inequality (up to localization) and the BDG inequality :
    \begin{align*}
        & \E\left[ \sup\limits_{t\leq T} \left\{ e^{\alpha t} |Y_t|^p + \frac{p(p-1)}{4} \int_t^T e^{\alpha s} |Y_s|^{p-2} \left(|Z_s|^2 + |Z_s^0|^2\right) \dd s \right\} \right] \\ 
        &\hspace*{1cm} \leq C_{p,T} \, \E\left[ |\xi|^p + \int_0^T |f(s,0,0,0,\delta_0)|^p \, \dd s \right] + \E \left[\sup\limits_{t\leq T} - p \int_t^T  e^{\alpha s} |Y_s|^{p-2} Y_s\cdot \left( Z_s\, \dd W_s + Z_s^0\, \dd W_s^0\right) \right]
        \\ &\hspace*{1cm} \leq C_{p,T} \, \E\left[ |\xi|^p + \int_0^T |f(s,0,0,0,\delta_0)|^p \,\dd s + \left( \int_0^T  e^{2\alpha s} |Y_s|^{2p-2} \left( |Z_s|^2 + |Z_s^0|^2\right) \dd s\right)^{1/2} \right]\\
        %&\hspace*{1cm} \leq C_{p,T}  \E\left[ |\xi|^p + \int_0^T |f(s,0,0,0,\delta_0)|^p \,\dd s \right] + \E \left[ \sup\limits_{s\leq T} \left\{ e^{\alpha s/2} |Y_s|^{p/2} \right\} \left( \int_0^T  e^{\alpha s} |Y_s|^{p-2} \left( |Z_s|^2 + |Z_s^0|^2\right) \dd s\right)^{1/2} \right]\\
         &\hspace*{1cm} \leq C_{p,T}\left\{ \E\left[ |\xi|^p + \int_0^T |f(s,0,0,0,\delta_0)|^p \,\dd s \right] \right. \\
         &\hspace*{4cm}+ \left. \E \left[ \varepsilon \sup\limits_{s\leq T} \left\{ e^{\alpha s} |Y_s|^{p} \right\} + \frac{1}{\varepsilon} \int_0^T  e^{\alpha s} |Y_s|^{p-2} \left( |Z_s|^2 + |Z_s^0|^2\right) \dd s \right] \right\}.
    \end{align*}
    where we handle the last term with Young's inequality and the previous inequality. We eventually get the expected result :
    \[\E\left[ \sup\limits_{t\leq T} \left\{ e^{\alpha t} |Y_t|^p \right\} \right] \leq C_{p,T}\, \E\left[ |\xi|^p + \int_0^T |f(s,0,0,0,\delta_0)|^p\, \dd s \right] < +\infty.\]
\end{proof}

\section{Convergence of empirical distributions}

\noindent We now take a deeper look at the corresponding interacting particle system. Our goal is to show that the empirical measure converges to the conditional law with respect to $W^0$ of the underlying problem as the size of the population goes to infinity. For that, we consider the following system, for $i\in\{1, \dots, n\}$ :
\begin{align}\label{eqn1}
    \begin{split}
        Y_t^{i,n} &= \xi^i + \int_t^T f(s,Y_s^{i,n},Z_s^{i,i,n},Z_s^{0,i,n},\mu_s^n)\, \dd s - \int_t^T Z_s^{i,n}\cdot \dd \mathbf{W}_s - \int_t^T Z_s^{0,i,n}\, \dd W_s^0 \\
        &= \xi^i + \int_t^T f(s,Y_s^{i,n},Z_s^{i,i,n},Z_s^{0,i,n},\mu_s^n)\, \dd s - \int_t^T \sum\limits_{k=1}^n Z_s^{i,k,n}\, \dd W_s^k - \int_t^T Z_s^{0,i,n}\, \dd W_s^0
    \end{split}
\end{align}
where the $\xi^i$ are i.i.d. terminal conditions, the $(W_k)_{1\leq k\leq n}$ and $W^0$ are $d$-dimensional independent Brownian motions on $\Omega^1$ and $\Omega^0$ respectively, and we used the notations :
\[\mu_s^n = \frac{1}{n}\sum\limits_{i=1}^n \delta_{Y_t^{i,n}},~ \mathbf{W} = (W_k)_{1\leq k\leq n},~ Z^{i,n} = (Z^{i,1,n},\dots, Z^{i,n,n}). \]
To simplify, let us write the system above as the following equation in $(\R^m)^n$ :
\[\mathbf{Y}_t = \Xi + \int_t^T \mathbf{F}_s(\mathbf{Y}_s, \mathbf{Z}_s^{},\mathbf{Z}_s^0)\, \dd s - \int_t^T \mathbf{Z}_s\,\dd \mathbf{W}_s - \int_t^T \mathbf{Z}_s^0\, \dd \mathbf{W}_s^0,\]
where we denote $\mathbf{W}^0 = (W^0,\dots, W^0)$, $\mathbf{Z} = (Z^{1,n},\dots, Z^{n,n})$, $\mathbf{Z}^0 = (Z^{0,1,n},\dots, Z^{0,n,n})$, $\Xi = (\xi^1,\dots,\xi^n)$. The function $\mathbf{F}$ is defined as follows :
\[\mathbf{F} : (t,y,z,\tilde{z})\in [0,T]\times (\R^m)^n \times (\R^{m\times d})^{n\times n} \times (\R^{m\times d})^n \mapsto (f(t, y^i, (z^{diag})^{i,i}, \tilde{z}^{i,i}, \mu_y^n))_{1\leq i\leq n} \]
with $\mu_y^n = n^{-1} \sum \delta_{y^i}$ the empirical measure associated with $y$. We also used the notation $\mathbf{z}^{diag}$ to signal that the diagonal terms are the only left in the argument of $\mathbf{F}$ :
\[\mathbf{z}^{diag} = (z^{1,1,n},\dots, z^{n,n,n}) \in (\R^{m\times d})^n.\]
Well-posedness as well as existence and uniqueness follow from the classical result applied to $\mathbf{F}$ since for all $t\in [0,T]$, $y_1,y_2\in (\R^m)^n$ and $z_1,z_2,\tilde{z}_1,\tilde{z}_2$ (to make the expressions lighter, we write $z^{i,i} = (z^{diag})^{i,i}$) :
\begin{align*}
    |\mathbf{F}(t,y_1,z_1,\tilde{z}_1) - \mathbf{F}(t,y_2,z_2,\tilde{z}_2)|^2 &= \sum\limits_{i=1}^n |f(t, y_1^i, z_1^{i,i}, \tilde{z}_1^{i,i}, \mu_{y_1}^n) - f(t, y_2^i, z_2^{i,i}, \tilde{z}_2^{i,i}, \mu_{y_2}^n)|^2 \\
    &\leq C_f^2 \sum\limits_{i=1}^n \left( |y_1^i - y_2^i| + |z_1^{i,i} -z_2^{i,i}| + |\tilde{z}_1^{i,i} -\tilde{z}_2^{i,i}| + W_2(\mu_{y_1}^n, \mu_{y_2}^n) \right)^2\\
    &\leq 8C_f^2 \sum\limits_{i=1}^n \left( |y_1^i - y_2^i|^2 + |z_1^{i,i} -z_2^{i,i}|^2 + |\tilde{z}_1^{i,i} -\tilde{z}_2^{i,i}|^2 + \frac{1}{n} \sum\limits_{k=1}^n |y_1^k-y_2^k|^2 \right) \\
    &\leq 16C_f^2 \left( |y_1 - y_2|^2 + |z_1 - z_2|^2 + |\tilde{z}_1 - \tilde{z}_2|^2 \right).
\end{align*}

\noindent We then turn to the moment properties of the interacting particles system. To do so, we introduce the coupled system, which will be useful in the following :
\begin{equation}\label{eqn3}
    \tilde{Y}_t^i = \xi^i + \int_t^T f(s, \tilde{Y}_s^i, \tilde{Z}_s^i, \tilde{Z}_s ^{0,i}, \LL^1 (Y_s))\, \dd s - \int_t^T \tilde{Z}_s^i\,\dd W_s^i - \int_t^T \tilde{Z}_s^{0,i}\,\dd W_s^0 , ~~i=1,\dots,N.
\end{equation}
where $Y$ is the first component of the solution to the corresponding equation \eqref{eq1} given by the first theorem. Existence and uniqueness follow from the first result. The main difference with \eqref{eqn1} is the fact that the $\tilde{Y}^i$ are now i.i.d. conditionnally to $\FF^0$, the filtration associated with the common noise. In the following proposition, we drop the superscript $i$ for the sake of simplicity.

\begin{Prop}[] \label{PropUtile}
    Assume that $(f(s,0,0,0,\delta_0))\in \HH^{2p}$ for $p\geq 2$ and that
    \[\xi\in L^{2p},~~ \sup\limits_{t\leq T} \E\left[ |\tilde{Z}_t|^{2p} + |\tilde{Z}_t^0|^{2p} \right] < +\infty .\]
    There exists $C_0 > 0$ such that :
    \[\begin{array}{rll}
        \mathbb{E}\left[\sup\limits_{t\leq T} \E^0\left[ |\tilde{Y}_t|^{2p}\right]%^{2p/q}
        \right] &\leq C_0 & \\[5pt]
        \E \left[ |\tilde{Y}_s - \tilde{Y}_r|^p |\tilde{Y}_t - \tilde{Y}_s|^p\right] &\leq C_0 |t-r|^2,~~ &\text{for}~ 0\leq r\leq s\leq t\leq T\\[5pt]
        \E \left[ |\tilde{Y}_t - \tilde{Y}_s|^p\right] &\leq C_0 |t-s|,~~ &\text{for}~ 0\leq s\leq t\leq T \\[5pt]
        \E\left[ \sup\limits_{s\leq u \leq t} \E^0\left[ |\tilde{Y}_u - \tilde{Y}_s|^2 \right]^{p} \right] &\leq C_0 |t-s|^{p/2},~~ &\text{for}~ 0\leq s\leq t\leq T.
    \end{array}\]
\end{Prop}

\begin{proof}
    The first inequality is a consequence of the arguments given in the proof of Proposition \ref{PropLp}. Indeed, with the assumptions,
    \[\mathbb{E}\left[\sup\limits_{t\leq T} \E^0\left[ |\tilde{Y}_t|^{2p}\right]\right] \leq \mathbb{E}\left[ \E^0\left[ \sup\limits_{t\leq T} |\tilde{Y}_t|^{2p}\right]\right] = \E \left[ \sup\limits_{t\leq T} |\tilde{Y}_t|^{2p} \right] < +\infty.\]
    Let us turn to the second inequality. For simplicity, let us now denote $f_u = f(u, \tilde{Y}_u^i,\tilde{Z}_u^i,\tilde{Z}_u ^{0,i}, \LL^1 (Y_u))$ and $f_u^0 = f(u,0,0,0,\delta_0)$ for $0\leq u\leq T$ in the following. Using Hölder and Jensen inequalities, we get :
    \begin{align*}
        \E\left[ |\tilde{Y}_s - \tilde{Y}_r|^p |\tilde{Y}_t - \tilde{Y}_s|^p\right] &\leq C_p \Bigg\{ (t-s)^{p-1} (s-r)^{p-1} \E\left[ \int_s^t |f_u|^p\, \dd u \int_r^s |f_u|^p\, \dd u \right] \\
        & \hspace*{-1cm}+ (t-s)^{p-1} \E\left[ \int_s^t |f_u|^p\, \dd u \left( \left| \int_r^s \tilde{Z}_u \, \dd W_u\right|^p + \left| \int_r^s \tilde{Z}_u^0 \, \dd W_u^0\right|^p \right) \right] \\
        & \hspace*{-1cm}+ (s-r)^{p-1} \E\left[ \int_r^s |f_u|^p\, \dd u \left( \left| \int_s^t \tilde{Z}_u \, \dd W_u\right|^p + \left| \int_s^t \tilde{Z}_u^0 \, \dd W_u^0\right|^p \right) \right] \\
        & \hspace*{-1cm}+ \E\left[ \left| \int_s^t \tilde{Z}_u \, \dd W_u\right|^{2p} \right]^{1/2} \left( \E\left[ \left| \int_r^s \tilde{Z}_u \, \dd W_u\right|^{2p} \right]^{1/2} + \E\left[ \left| \int_r^s \tilde{Z}_u^0 \, \dd W_u^0\right|^{2p} \right]^{1/2}\right)\\
        & \hspace*{-1cm} + \E\left[ \left| \int_s^t \tilde{Z}_u^0 \, \dd W_u^0\right|^{2p} \right]^{1/2} \left( \E\left[ \left| \int_r^s \tilde{Z}_u \, \dd W_u\right|^{2p} \right]^{1/2} + \E\left[ \left| \int_r^s \tilde{Z}_u^0 \, \dd W_u^0\right|^{2p} \right]^{1/2}\right) \Bigg\}.
    \end{align*}
    Now thanks to the Burkholder-Davis-Gundy inequality, the above becomes :
    \begin{align*}
        \E\left[ |\tilde{Y}_s - \tilde{Y}_r|^p |\tilde{Y}_t - \tilde{Y}_s|^p\right] &\leq C_p \Bigg\{ (t-s)^{p-1} (s-r)^{p-1} T \E\left[ \left( \int_0^T |f_u|^{2p}\, \dd u\right) \right] \\
        & \hspace*{-1cm}+ (t-s)^{p-1/2} \E\left[ \int_s^t |f_u|^{2p}\, \dd u\right]^{1/2} \left( \E\left[ \left| \int_r^s \tilde{Z}_u \, \dd W_u\right|^{2p}\right]^{1/2} + \E\left[ \left| \int_r^s \tilde{Z}_u^0 \, \dd W_u^0\right|^{2p}\right]^{1/2} \right) \\
        & \hspace*{-1cm}+ (s-r)^{p-1/2} \E\left[ \int_r^s |f_u|^{2p}\, \dd u\right]^{1/2} \left( \E\left[ \left| \int_s^t \tilde{Z}_u \, \dd W_u\right|^{2p}\right]^{1/2} + \E\left[ \left| \int_s^t \tilde{Z}_u^0 \, \dd W_u^0\right|^{2p}\right]^{1/2} \right) \\
        & \hspace*{-1cm}+ \E\left[ \left(\int_s^t |\tilde{Z}_u|^2 \, \dd u\right)^{p} \right]^{1/2} \left( \E\left[ \left( \int_r^s |\tilde{Z}_u|^2 \, \dd u\right)^{p} \right]^{1/2} + \E\left[ \left( \int_r^s |\tilde{Z}_u^0|^2 \, \dd u\right)^{p} \right]^{1/2}\right)\\
        & \hspace*{-1cm}+ \E\left[ \left(\int_s^t |\tilde{Z}_u^0|^2 \, \dd u\right)^{p} \right]^{1/2} \left( \E\left[ \left( \int_r^s |\tilde{Z}_u|^2 \, \dd u\right)^{p} \right]^{1/2} + \E\left[ \left( \int_r^s |\tilde{Z}_u^0|^2 \, \dd u\right)^{p} \right]^{1/2}\right) \Bigg\}.
    \end{align*}
    We use the following, coming from the Lipschitz hypothesis and the additional assumptions :
    \begin{align*}
        \E \left[ \int_0^T |f_u|^{2p}\, \dd u\right] &\leq C_p \E\left[ \int_0^T \left( |f_u^0|^{2p} + |\tilde{Y}_u|^{2p} + |\tilde{Z}_u|^{2p} + |\tilde{Z}_u^0|^{2p} + W_2\left(\LL^1(Y_u),\delta_0\right)^{2p} \right) \, \dd u \right]\\
        &\leq C_{p,T} \left( \E \left[ \int_0^T |f_u^0|^{2p}\, \dd u\right] + \sup\limits_{t\leq T} \E\left[ |\tilde{Y}_t|^{2p}\right] + \sup\limits_{t\leq T} \E\left[ |Y_t|^{2p}\right] \right. \\
        &\hspace*{3cm} \left. \phantom{\E \left[ \int_0^T \right]}+ \sup\limits_{t\leq T} \E\left[ |\tilde{Z}_t|^{2p}\right] + \sup\limits_{t\leq T} \E\left[ |\tilde{Z}_t^0|^{2p}\right]  \right) < +\infty
    \end{align*}
    as well as :
    \[ \E\left[ \left(\int_s^t |\tilde{Z}_u|^2 \, \dd u\right)^{p} \right] \leq (t-s)^{p-1} \E \left[ \int_s^t |\tilde{Z}_u|^{2p}\, \dd u\right] \leq C_{Z} (t-s)^{p}.\]
    Coming back to the desired inequality :
    \begin{align*}
        \E\left[ |\tilde{Y}_s - \tilde{Y}_r|^p |\tilde{Y}_t - \tilde{Y}_s|^p\right] &\leq \tilde{C}_p \Bigg\{ \E\left[ (t-s)^{p-1} (s-r)^{p-1} T \left( \int_0^T |f_u|^{2p}\, \dd u\right) \right] \\
        & \hspace*{1cm}+ (t-s)^{p-1/2}\, \E\left[ \int_s^t |f_u|^{2p}\, \dd u\right]^{1/2} (s-r)^{\frac{p}{2}}  \\
        & \hspace*{1cm}+ (s-r)^{p-1/2}\, \E\left[ \int_r^s |f_u|^{2p}\, \dd u\right]^{1/2} (t-s)^{\frac{p}{2}} \\
        & \hspace*{1cm}+ (t-s)^{\frac{p}{2}} (s-r)^{\frac{p}{2}} + (t-s)^{\frac{p}{2}} (s-r)^{\frac{p}{2}} \Bigg\}.\\
        &\leq C_{p,T} (t-r)^{p} \leq C_{p,T} (t-r)^{2}
    \end{align*}
    since $p\geq 2$. Similarly, for the third one, we have :
    \begin{align*}
        \E\left[ |\tilde{Y}_t - \tilde{Y}_s|^p\right] &\leq C_p \left\{ \E\left[ \int_s^t (t-s)^{p-1}| f(u,0,0,0, \delta_0)|^p \, \dd u \right] \phantom{+ \E \left[ \left( \int_s^t \right) ^{p/2} \right]}\right. \\
        & \hspace*{1cm}+ (t-s) \E\left[ \int_s^t |\tilde{Y}_u|^p + |\tilde{Z}_u|^p + |\tilde{Z}_u^0|^p + W_2\left(\LL^1(Y_u,\delta_0)\right)^p \dd u \right] \\
        &\left. \hspace*{2cm}+ \E \left[ \left( \int_s^t |\tilde{Z}_u|^2\, \dd u\right) ^{p/2} \right] + \E \left[ \left( \int_s^t |\tilde{Z}_u^0|^2\, \dd u\right) ^{p/2} \right] \right\}\\
        &\leq C_{p,T} \left\{ (t-s)^{p-1} + (t-s) + (t-s)^{p/2} \right\} \leq C_{p,T} (t-s).
    \end{align*}
    For the last inequality, we use similar arguments to get the following :
    \begin{align*}
        \E\left[ \sup\limits_{s\leq u \leq t} \E^0\left[ |\tilde{Y}_u - \tilde{Y}_s|^2 \right]^{p} \right] &\leq C_p \E\left[ \sup\limits_{s\leq u \leq t} \E^0\left[ \left| \int_s^u f_v \, \dd v\right|^{2p} \right] + \E^0 \left[ \sup\limits_{s\leq u \leq t} \left|\int_u^s \tilde{Z}_s \, \dd W_s\right|^{2p} \right] \right. \\
        &\hspace*{4cm} \left. + \E^0 \left[ \sup\limits_{s\leq u \leq t} \left|\int_u^s \tilde{Z}_s^0 \, \dd W_s^0\right|^{2p} \right] \right]\\
        &\leq C_p \left\{ (t-s)^{2p-1} \E\left[ \int_s^t |f_v|^{2p}\, \dd p \right] + (t-s)^p \right\} \leq C_0 |t-s|^{p/2}
    \end{align*}
    using the Burkholder-Davis-Gundy inequality, as well as the control on terms already given above. This concludes the proof.
\end{proof}

\noindent In the following, we sometimes make use of this additional hypothesis on the $Z$-component of the solution.
\begin{enumerate}
    \item[\textbf{(H3)}] There exists $p\geq 2$ such that the driver satisfies $(f(s,0,0,0,\delta_0))\in \HH^{2p}$ and that
    \[\xi\in L^{2p},~~ \sup\limits_{t\leq T} \E\left[ |\tilde{Z}_t|^{2p} + |\tilde{Z}_t^0|^{2p} \right] < +\infty .\]
\end{enumerate}

\begin{Rmq}
    As underlined in \cite{LauriereTangpi}, the assumption on the moments of $(Z_t)$ and $(Z_t^0)$ holds in usual cases, described in \cite{ImkellerDosReis} and \cite{CheriditoNam}.
\end{Rmq}

\noindent Adapting the proof of Lemma 17 in \cite{LauriereTangpi}, we get the following control for the empirical measure of the system \eqref{eqn1}.

\begin{Lem}\label{Interm1}
    Let us denote $(Y^{i,n},Z^{i,n},Z^{0,i,n})_{1\leq i\leq n}$ (resp. $(\tilde{Y}^{i},\tilde{Z}^{i},\tilde{Z}^{0,i})_{1\leq i\leq n}$) a system of solutions to \eqref{eqn1} (resp. \eqref{eqn3}) as well as $(Y,Z,Z^0)$ a solution to \eqref{eq1}. Under the standard assumptions (H1)-(H2),
    \begin{equation} \label{eqnG}
        W_p \left( \frac{1}{n} \sum\limits_{i=1}^n \delta_{Y_t^{i,n}}, \LL^1(Y_t) \right) \leq \exp\left( T e^{C_f T}\right) W_p \left( \frac{1}{n} \sum\limits_{i=1}^n \delta_{\tilde{Y}_t^i}, \LL^1(Y_t) \right).
    \end{equation}
\end{Lem}

\begin{proof}
    Coming back to \eqref{eqn1} and \eqref{eqn3}, we first write, for $i\in\{1,\dots,n\}$ fixed :
    \begin{align*}
        \tilde{Y}_t^i - Y_t^{i,n} &= \int_t^T f(s, \tilde{Y}_s^i, \tilde{Z}_s^i, \tilde{Z}_s ^{0,i}, \LL^1 (Y_s)) - f(s,Y_s^{i,n},Z_s^{i,i,n},Z_s^{0,i,n},\mu_s^n)\, \dd s \\
        &\hspace*{3cm}- \sum\limits_{k=1}^n \int_t^T \left( \tilde{Z}_s^i \delta_{ik} - Z_s^{i,k,n} \right) \dd W_s^k - \int_t^T (\tilde{Z}_s^{0,i} - Z_s^{0,i,n}) \,\dd W_s^0 \\
        %&= \int_t^T f(s, \tilde{Y}_s^i, \tilde{Z}_s^i, \tilde{Z}_s ^{0,i}, \LL^1 (Y_s)) - f(s, Y_s^{i,n}, \tilde{Z}_s^i, \tilde{Z}_s ^{0,i}, \LL^1 (Y_s))\\
        %&\hspace*{1cm} + f(s, Y^{i,n}, \tilde{Z}_s^i, \tilde{Z}_s ^{0,i}, \LL^1 (Y_s)) - f(s, Y_s^{i,n}, Z_s^{i,i,n}, \tilde{Z}_s ^{0,i}, \LL^1 (Y_s)) \\
        %&\hspace*{1cm} + f(s, Y_s^{i,n}, Z_s^{i,i,n}, \tilde{Z}_s ^{0,i}, \LL^1 (Y_s)) - f(s, Y_s^{i,n}, Z_s^{i,i,n}, Z_s^{0,i,n}, \LL^1 (Y_s)) \\
        %&\hspace*{1cm}+ f(s, Y_s^{i,n}, Z_s^{i,i,n}, Z_s^{0,i,n}, \LL^1 (Y_s)) - f(s,Y_s^{i,n},Z_s^{i,i,n}, Z_s^{0,i,n}, \mu_s^n)\, \dd s \\
        %&\hspace*{3cm}- \sum\limits_{k=1}^n \int_t^T \left( \tilde{Z}_s^i \delta_{ik} - Z_s^{i,k,n} \right) \dd W_s^k - \int_t^T \tilde{Z}_s^{0,i} - Z_s^{0,i,n} \,\dd W_s^0 \\
        &= \int_t^T \alpha_s^i + \beta_s^i + \gamma_s^i (\tilde{Z}_s^i - Z_s^{i,i,n}) + \gamma_s^{0,i} (\tilde{Z}_s^{0,i} - Z_s^{0,i,n}) \, \dd s \\
        &\hspace*{3cm} - \sum\limits_{k=1}^n \int_t^T \left( \tilde{Z}_s^i \delta_{ik} - Z_s^{i,k,n} \right) \dd W_s^k - \int_t^T (\tilde{Z}_s^{0,i} - Z_s^{0,i,n}) \,\dd W_s^0
    \end{align*}
    where we used the following, given by a Taylor formula :
    \begin{align*}
        \alpha_s^i &= f(s, \tilde{Y}_s^i, \tilde{Z}_s^i, \tilde{Z}_s ^{0,i}, \LL^1 (Y_s)) - f(s, Y_s^{i,n}, \tilde{Z}_s^i, \tilde{Z}_s ^{0,i}, \LL^1 (Y_s))\\
        \beta_s^i &= f(s, Y_s^{i,n}, Z_s^{i,i,n}, Z_s^{0,i,n}, \LL^1 (Y_s)) - f(s,Y_s^{i,n},Z_s^{i,i,n}, Z_s^{0,i,n}, \mu_s^n)\\
        \gamma_s^i &= \int_0^1 \partial_z f (s, Y_s^{i,n}, Z_s^{i,i,n} + \lambda (\tilde{Z}_s^i - Z_s^{i,i,n}), \tilde{Z}_s^{0,i}, \LL^1(Y_s))\, \dd \lambda \\
        \gamma_s^{0,i} &= \int_0^1 \partial_{\tilde{z}} f (s, Y_s^{i,n}, Z_s^{i,i,n}, Z_s^{0,i,n} + \lambda (\tilde{Z}_s^{0,i} - Z_s^{0,i,n}), \LL^1(Y_s))\, \dd \lambda .
    \end{align*}
    Since $f$ is Lipschitz, the stochastic exponential $\EE(\gamma^i \cdot W^i + \gamma^{0,i}\cdot W^0)$ defines an equivalent probability measure~$\Q$. The square integrability of the variables $\tilde{Z}^i, \tilde{Z}^{0,i}, Z^{0,i,n}$ and $Z^{i,k,n}$ ensure that, taking the expectation with respect to $\Q$ thanks to Girsanov's theorem :
    \begin{align*}|\tilde{Y}_t^i - Y_t^{i,n}| &= \left| \E_\Q \left[ \int_t^T \alpha_s^i + \beta_s^i \, \dd s \mid \FF_t^n\right]\right|\\
    & \leq C_f \E_\Q \left[ \int_t^T \left| \tilde{Y}_s^i - Y_s^{i,n}\right| + W_2(\LL^1(Y_s),\mu_s^n)\, \dd s\mid \FF_t^n\right].
    \end{align*}
    using the Lipschitz property. Now by Gronwall's inequality, it follows that :
    \[|\tilde{Y}_t^i - Y_t^{i,n}|\leq e^{C_f T} \E_\Q \left[ \int_t^T W_2 (\LL^1(Y_s), \mu_s^n)\, \dd s\mid \FF_t^n\right]. \]
    On the one hand, the triangle inequality implies :
    \[W_p \left( \frac{1}{n} \sum\limits_{i=1}^n \delta_{Y_t^{i,n}}, \LL^1(Y_t) \right) \leq W_p \left( \frac{1}{n} \sum\limits_{i=1}^n \delta_{Y_t^{i,n}}, \frac{1}{n} \sum\limits_{i=1}^n \delta_{\tilde{Y}_t^i} \right) + W_p \left( \frac{1}{n} \sum\limits_{i=1}^n \delta_{\tilde{Y}_t^i}, \LL^1(Y_t) \right), \]
    and on the other hand, with the previous inequalities, since $W_2 \leq W_p$ :
    \[W_p \left( \frac{1}{n} \sum\limits_{i=1}^n \delta_{Y_t^{i,n}}, \frac{1}{n} \sum\limits_{i=1}^n \delta_{\tilde{Y}_t^i} \right) \leq \left( \frac{1}{n} \sum\limits_{i=1}^n |\tilde{Y}_t^i - Y_t^{i,n}|^p \right)^{1/p} \leq e^{C_f T} \E_\Q \left[ \int_t^T W_2 (\LL^1(Y_s), \mu_s^n)\, \dd s\mid \FF_t^n\right].\]
    Coming back to the triangle inequality and injecting this last inequality, we can now apply Gronwall's lemma to get \eqref{eqnG}.
\end{proof}

\begin{Prop}[A first result of conditional propagation of chaos]\label{1PoC}
    For $p >4$, assume furthermore that $\xi\in L^p$ and that $(f(s,0,0,0,\delta_0))\in \HH^p$. Then there exists a constant $C = C(d,p,T)$ such that :
    \[\E \left[ \sup\limits_{t\leq T} \E^0\left[ W_2^2(\mu_t^n,\mu_t) \right] \right] \leq C\, \varepsilon_n = C \times  \begin{cases} n^{-1/2} & \text{ if } d < 4, \\ n^{-1/2} \log(n) & \text{ if } d = 4, \\ n^{-2/d} & \text{ if } d > 4. \end{cases}\]
\end{Prop}

\begin{proof}
    We derive from \cite[Theorem 5.8]{CD2} and \cite[Remark 5.9]{CD2} the following inequality :
    \begin{equation}
        \E^0 \left[ W_2^2 (\mu_t^n,\mu_t) \right] \leq C(d,p) M_q(\mu_t)^2 \varepsilon_n.
    \end{equation}
    Taking the supremum in time and the expectation with respect to the whole probability space, we get the desired result thanks to Proposition \ref{PropLp}.
\end{proof}

\noindent Under the more restrictive assumptions of Proposition \ref{PropUtile}, the Lemma 1 from \cite{BCCdRH}, based on a result from \cite{FournierGuillin}, holds and yields the following result.

\begin{Thm}\label{ThmChaos}
    Under (H1)-(H2)-(H3), recalling that $\mu_s = \LL^1(Y_s)$ and $\mu_s^n = n^{-1} \sum \delta_{Y_s^{i,n}}$, there exists a constant $C > 0$ such that :
    \begin{equation}
        \E \left[ \sup\limits_{s\leq T} W_2^2 (\mu_s^n, \mu_s) \right] \leq C\tilde{\varepsilon}_n = C \times  \begin{cases} n^{-1/2 + 1/p} & \text{ if } d < 4, \\ n^{-1/2 + 1/p} \log(1+n)^{1-2/p} & \text{ if } d = 4, \\ n^{-2(1-2/p)/d} & \text{ if } d > 4. \end{cases}
    \end{equation}
\end{Thm}

\begin{proof}
    The result is an immediate consequence of the Lemma 1 from \cite{BCCdRH}, given by Proposition \ref{PropUtile} combined with Lemma \ref{Interm1}.
\end{proof}

\section{Interacting particles approximation}
\begin{Lem}
    Under the assumptions (H1)-(H2), we have the following moment estimates :
    \begin{equation}
        \E \left[ \frac{1}{n} \sum\limits_{i=1}^n |Y_t^{i,n}|^2 + \frac{1}{n} \sum\limits_{i,k=1}^n \int_0^T |Z_s^{i,k,n}|^2\, \dd s + \frac{1}{n} \sum\limits_{i=1}^n \int_0^T |Z_s^{0,i,n}|^2\, \dd s \right] \leq C_T.
    \end{equation}
\end{Lem}

\begin{proof}
    Applying Itô's formula to the family $e^{\alpha t} |Y_t^{i,n}|^2$, we get after summing over $i\in\{1,\dots, n\}$ :
    \begin{align*}
    \frac{1}{n} &\sum_{i=1}^n e^{\alpha t} |Y_t^{i,n}|^2 + \frac{1}{n} \sum\limits_{i,k=1}^n \int_t^T |Z_s^{i,k,n}|^2\, \dd s + \frac{1}{n} \sum\limits_{i=1}^n \int_t^T |Z_s^{0,i,n}|^2\, \dd s \\
    &= \frac{1}{n} \sum_{i=1}^n e^{\alpha t} |\xi^i|^2 - \frac{\alpha}{n} \sum_{i=1}^n \int_t^T e^{\alpha s} |Y_s^{i,n}|^2 \, \dd s + \frac{2}{n} \sum_{i=1}^n \int_t^T e^{\alpha s} Y_s^{i,n} \cdot f(s,Y_s^{i,n},Z_s^{i,i,n},Z_s^{0,i,n},\mu_s^n)\, \dd s\\
    &- \frac{2}{n} \sum_{i=1}^n \int_t^T e^{\alpha s} Y_s^{i,n} \cdot \left( \sum_{k=1}^n Z_s^{i,k,n}\, \dd W_s^k + Z_s^{0,i,n}\, \dd W_s^0\right).
    \end{align*}
    Using the Lipschitz property on $f$ as well as Young's inequality gives :
    \begin{align*}
        &\int_t^T e^{\alpha s} Y_s^{i,n} \cdot f(s,Y_s^{i,n},Z_s^{i,i,n},Z_s^{0,i,n},\mu_s^n)\, \dd s \\
        &\leq C\left[ \int_t^T e^{\alpha_s} |f_s^0|^2\, \dd s + \int_t^T e^{\alpha s} \left\{\left( 1 + \frac{1}{\varepsilon} \right)  |Y_s^{i,n}|^2 + \varepsilon \left( |Z_s^{i,i,n}|^2 + |Z_s^{0,i,n}|^2 + W_2^2 (\mu_s^n,\delta_0)\right) \right\} \dd s \right]
    \end{align*}
    where we denote $f_s^0 = f(s,0,0,0,\delta_0)$. Summing over $i$, choosing $\varepsilon$ small enough, then $\alpha$ big enough and taking expectation on both sides of the inequality yields the expected result, since :
    \[\E \left[ W_2^2 (\mu_s^n,\delta_0)\right] \leq \E\left[  \frac{1}{n} \sum_{i=1}^n |Y_s^{i,n}|^2 \right]. \]
\end{proof}

\begin{Thm}
    Under the same assumptions (H1)-(H2), there exists $C_T > 0$ such that :
    \begin{align*}
        &\E\left[ \sup\limits_{t\leq T} \frac{1}{n} \sum \limits_{i=1}^n \left( |Y_t^{i,n} - \tilde{Y}_t^i|^2 + \int_t^T |Z_s^{0,i,n} - \tilde{Z}_s^{0,i} |^2\, \dd s + \int_t^T \sum\limits_{k=1}^n | Z_s^{i,k,n} - \tilde{Z}_s^i \delta_{i,k} |^2\, \dd s \right) \right] \\ &\hspace*{4cm} \leq C_T \E \left[ \sup\limits_{s\leq T} W_2^2 (\mu_s^n, \mu_s) \right],\\
        &\E\left[ \sup\limits_{t\leq T}  \left\{ |Y_t^{i,n} - \tilde{Y}_t^i|^2 + \int_t^T |Z_s^{0,i,n} - \tilde{Z}_s^{0,i} |^2\, \dd s + \int_t^T \sum\limits_{k=1}^n | Z_s^{i,k,n} - \tilde{Z}_s^i \delta_{i,k} |^2\, \dd s \right\} \right] \\ &\hspace*{4cm} \leq C_T \E \left[ \sup\limits_{s\leq T} W_2^2 (\mu_s^n, \mu_s) \right].
    \end{align*}
\end{Thm}

\begin{proof}
    We first derive from Ito's formula the following identity for $i\in\{ 1,\dots,n\}$:
    \begin{align*} 
        &e^{\alpha t} |Y_t^{i,n} - \tilde{Y}_t^i |^2 + \int_t^T e^{\alpha s} \left( \sum\limits_{k=1}^n |Z_s^{i,k,n} - \tilde{Z}_s^i | ^2 + |Z_s^{0,i,n} - \tilde{Z}_s^{0,i} |^2 \right)\, \dd s \\  &= -\alpha \int_t^T e^{\alpha s} |Y_s^{i,n} - \tilde{Y}_s^i |^2 \, \dd s - 2 \int_t^T e^{\alpha s} (Y_s^{i,n} - \tilde{Y}_s^i)\cdot \dd (Y_s^{i,n} - \tilde{Y}_s^i). 
    \end{align*}
    Repeating the arguments of the previous proof, using the Lipschitz continuity of $f$  and Young's inequality we get for $\alpha$ large enough :
    \begin{align} \label{Ineq1} \small
    \begin{split}
        &e^{\alpha t} |Y_t^{i,n}-\tilde{Y}_t^i|^2 + \sum\limits_{k=1}^n \int_t^T e^{\alpha s} |Z_s^{i,k,n} - \tilde{Z}_s^i \delta_{i,k}|^2\, \dd s + \int_t^T e^{\alpha s} |Z_s^{0,i,n}- \tilde{Z}_s^{0,i}|^2\, \dd s\\ & \leq C W_2^2(\mu_s^n,\mu_s) - 2 \int_t^T e^{\alpha s} (Y_s^{i,n} - \tilde{Y}_s^i) \sum\limits_{k=1}^n (Z_s^{i,k,n} - \tilde{Z}_s^{i} \delta_{i,k}) \, \dd W_s^k - 2 \int_t^T e^{\alpha s} (Y_s^{i,n} - \tilde{Y}_s^i) (Z_s^{0,i,n} - \tilde{Z}_s^{0,i}) \, \dd W_s^0,
        \end{split}
    \end{align}
    the two last terms of the above inequality being martingales. Summing over $i$ leads to the following :
    \begin{align}\label{IneqInterm2}
    \begin{split}
        \frac{1}{n} &\sum_{i=1}^n e^{\alpha t} |Y_t^{i,n}-\tilde{Y}_t^i|^2 + \frac{1}{n} \sum\limits_{i,k=1}^n \int_t^T e^{\alpha s} |Z_s^{i,k,n} - \tilde{Z}_s^i \delta_{i,k}|^2\, \dd s + \frac{1}{n} \sum\limits_{i=1}^n \int_t^T e^{\alpha s} |Z_s^{0,i,n}- \tilde{Z}_s^{0,i}|^2\, \dd s \\
        &\hspace*{9cm}\leq C W_2^2(\mu_t^n,\mu_t) + M_n(t,T) .
    \end{split}
    \end{align}
    where we denote by $M_n(t,T)$ the following martingale (on the product space $\Omega$) :
    \[M_n(t,T) = -\frac{2}{n} \sum\limits_{i=1}^n \int_t^T (Y_s^{i,n} - \tilde{Y}_s^i) \sum\limits_{k=1}^n (Z_s^{i,k,n} - \tilde{Z}_s^{i} \delta_{i,k}) \, \dd W_s^k -\frac{2}{n} \sum\limits_{i=1}^n \int_t^T (Y_s^{i,n} - \tilde{Y}_s^i) (Z_s^{0,i,n} - \tilde{Z}_s^{0,i}) \, \dd W_s^0\]
    Remark first that, taking the expectation, we get :
    \begin{equation}
        \begin{split}\label{IneqInterm1}
            &\E\left[ \frac{1}{n} \sum_{i=1}^n |Y_t^{i,n}-\tilde{Y}_t^i|^2 + \frac{1}{n} \sum\limits_{i,k=1}^n \int_t^T |Z_s^{i,k,n} - \tilde{Z}_s^i \delta_{i,k}|^2\, \dd s + \frac{1}{n} \sum\limits_{i=1}^n \int_t^T |Z_s^{0,i,n}- \tilde{Z}_s^{0,i}|^2\, \dd s \right] \\ & \hspace*{11cm}\leq C \E\left[ W_2^2(\mu_t^n,\mu_t) \right].
        \end{split}
    \end{equation}
    Note then that, thanks to the Burkhölder-Davis-Gundy inequality combined with the Cauchy-Schwarz inequality :
    \begin{align*}
        \E [ \sup\limits_{t\leq T} M_n(t,T)]
        %&\leq C \E \left[ \left( \frac{1}{n^2} \sum\limits_{j=1}^n \int_0^T \left( \sum\limits_{i=1}^n e^{\alpha s} |Y_s^{i,n} - \tilde{Y}_s^i| |Z_s^{i,k,n} - \tilde{Z}_s^{i} \delta_{i,k}|\right)^2 \dd s\right] \\  &\hspace*{5cm} +C\E \left[ + \frac{1}{n^2} \int_0^T \left( \sum\limits_{i=1}^n |Z_s^{0,i,n} - \tilde{Z}_s^{0,i}|^2 \right) \dd s \right)^{1/2} \right] \\
        &\leq C\E \left[ \left( \sum\limits_{k=1}^n \int_0^T \left( \frac{1}{n} \sum\limits_{i=1}^n e^{\alpha s} |Y_s^{i,n} - \tilde{Y}_s^i|^2 \right) \left( \frac{1}{n} \sum\limits_{i=1}^n e^{\alpha s} |Z_s^{i,k,n} - \tilde{Z}_s^{i} \delta_{i,k}|^2 \right) \dd s \right)^{1/2}\right] \\ & \hspace*{1.8cm} + C\E \left[ \left( \int_0^T \left( \frac{1}{n} \sum\limits_{i=1}^n e^{\alpha s} |Y_s^{i,n} - \tilde{Y}_s^i|^2 \right) \left( \frac{1}{n} \sum\limits_{i=1}^n e^{\alpha s} |Z_s^{0,i,n} - \tilde{Z}_s^{0,i}|^2 \right) \dd s \right)^{1/2}\right] \\
        &\leq 2C\varepsilon \E \left[ \sup\limits_{t\leq T} \left\{ \frac{1}{n} \sum\limits_{i=1}^n e^{\alpha s} |Y_t^{i,n} - \tilde{Y}_t^i|^2 \right\} \right] \\ & \hspace*{1.8cm} + \frac{C}{\varepsilon} \E\left[ \frac{1}{n} \sum\limits_{i=1}^n \int_0^T e^{\alpha s} \left( \sum\limits_{k=1}^n |Z_s^{i,k,n} - \tilde{Z}_s^{i} \delta_{i,k}|^2 + |Z_s^{0,i,n} - \tilde{Z}_s^{0,i}|^2 \right) \dd s \right].
    \end{align*}
    where we used Young's inequality for the last line. Let us now come back to \eqref{IneqInterm2} and recall \eqref{IneqInterm1}. Taking first the supremum and then the expectation, as well as using the previous inequality with $\varepsilon$ small enough (up to localization),
    \begin{align*}
        &\E\left[ \sup\limits_{t\leq T} \frac{1}{n} \sum \limits_{i=1}^n \left( |Y_t^{i,n} - \tilde{Y}_t^i|^2 + \int_t^T |Z_s^{0,i,n} - \tilde{Z}_s^{0,i} |^2\, \dd s + \int_t^T \sum\limits_{k=1}^n | Z_s^{i,k,n} - \tilde{Z}_s^i \delta_{i,k} |^2\, \dd s \right) \right] \\ &\hspace*{4cm} \leq C_T \E \left[ \sup\limits_{s\leq T} W_2^2 (\mu_s^n, \mu_s) \right].
    \end{align*}
    A similar reasoning starting from \eqref{Ineq1} instead gives the second desired inequality.
\end{proof}

\noindent Combined with Theorem \ref{ThmChaos}, this gives the following uniform convergence for the system of interacting particles.

\begin{Coro}\label{Corolla1}
    Under the stronger assumptions (H1)-(H2)-(H3) already used in Theorem 9, the following holds for $i\in\{1,\dots,n\}$ :
    \begin{align*}
        &\E\left[ \sup\limits_{t\leq T}  \left\{ |Y_t^{i,n} - \tilde{Y}_t^i|^2 + \int_t^T |Z_s^{0,i,n} - \tilde{Z}_s^{0,i} |^2\, \dd s + \int_t^T \sum\limits_{j=1}^n | Z_s^{i,k,n} - \tilde{Z}_s^i |^2\, \dd s \right\} \right] \\ &\hspace*{4cm} \leq C_T \E \left[ \sup\limits_{s\leq T} W_2^2 (\mu_s^n, \mu_s) \right] \leq C_T \tilde\varepsilon_n.
    \end{align*}
\end{Coro}

\noindent Mimicking the previous proof and using the same arguments as well as the independance between the $W^k$ and $W^0$, we also get the result with the supremum outside the $\E^0$ expectation, getting rid of the hypothesis (H3) in exchange for a less restrictive integrability condition.

\begin{Prop}\label{Propo13}
    If (H1) and (H2) hold, assuming also that $\xi\in L^p$ and that $(f(s,0,0,0,\delta_0))\in \HH^p$, for $p>4$, there exists $C_T > 0$ such that the following holds :
    \begin{align*}
        &\E\left[ \sup\limits_{t\leq T} \E^0 \left[\frac{1}{n} \sum \limits_{i=1}^n \left( |Y_t^{i,n} - \tilde{Y}_t^i|^2 + \int_t^T |Z_s^{0,i,n} - \tilde{Z}_s^{0,i} |^2\, \dd s + \int_t^T \sum\limits_{k=1}^n | Z_s^{i,k,n} - \tilde{Z}_s^i \delta_{i,k} |^2\, \dd s \right) \right] \right] \\ &\hspace*{9cm} \leq C_T \E \left[ \sup\limits_{s\leq T} \E^0 \left[ W_2^2 (\mu_s^n, \mu_s) \right] \right],\\
        &\E\left[ \sup\limits_{t\leq T} \E^0 \left[ \left\{ |Y_t^{i,n} - \tilde{Y}_t^i|^2 + \int_t^T |Z_s^{0,i,n} - \tilde{Z}_s^{0,i} |^2\, \dd s + \int_t^T \sum\limits_{k=1}^n | Z_s^{i,k,n} - \tilde{Z}_s^i \delta_{i,k} |^2\, \dd s \right\} \right] \right] \\ &\hspace*{9cm} \leq C_T \E \left[ \sup\limits_{s\leq T} \E^0 \left[ W_2^2 (\mu_s^n, \mu_s) \right] \right].
    \end{align*}
\end{Prop}

\noindent Recalling the Proposition \ref{1PoC}, we get an equivalent result for Corollary \ref{Corolla1} when the supremum and the expectation $\E^0$ are exchanged.

\begin{Coro}
    Under the same assumptions as Proposition \ref{Propo13}, the following holds for $i\in\{1,\dots,n\}$ :
    \begin{align*}
        &\E\left[ \sup\limits_{t\leq T}  \E^0 \left[ |Y_t^{i,n} - \tilde{Y}_t^i|^2 + \int_t^T |Z_s^{0,i,n} - \tilde{Z}_s^{0,i} |^2\, \dd s + \int_t^T \sum\limits_{j=1}^n | Z_s^{i,k,n} - \tilde{Z}_s^i |^2\, \dd s \right] \right] \\ &\hspace*{4cm} \leq C_T \E \left[ \sup\limits_{s\leq T} \E^0 \left[ W_2^2 (\mu_s^n, \mu_s) \right] \right] \leq C_T \varepsilon_n.
    \end{align*}
\end{Coro}

\bibliographystyle{plain}
\bibliography{biblio}

\end{document}